\input amstex
\documentstyle{amsppt}
\magnification1200
\tolerance=1000
\def\n#1{\Bbb #1}

\def\p{\Bbb C_{\infty}}

\def\lim{\hbox{lim }}

\def\al{\alpha}
\def\th{\theta}
\def\om{\omega}
\def\g{\goth }

\def\Hom{\hbox{Hom }}
\def\End{\hbox{End }}
\def\Ker{\hbox{Ker }}

\def\invlim{\hbox{invlim }}
\def\Spec{\hbox{Spec }}
\def\dim{\hbox{dim }}

\def\ord{\hbox{ord }}

\def\Corr{\hbox{ Corr}}
\def\Lie{\hbox{Lie}}
\def\diag{\hbox{ diag }}

\def\vf{\varphi}

\topmatter
\title
Reductions of Hecke correspondences on Anderson modular objects
\endtitle
\author
A. Grishkov, D. Logachev\footnotemark \footnotetext{E-mails: shuragri{\@}gmail.com; logachev94{\@}gmail.com (corresponding author)\phantom{***************}}
\endauthor
\NoRunningHeads
\thanks Thanks: The authors are grateful to FAPESP, S\~ao Paulo, Brazil for a financial support (process No. 2017/19777-6). The first author is grateful to SNPq, Brazil, to RFBR, Russia, grant 16-01-00577a (Secs. 1-4), and to Russian Science Foundation, project 16-11-10002 (Secs. 5-8) for a financial support. The second author is grateful to Laurent Lafforgue for invitation to
IHES where this paper was started. Discussions with Laurent Fargues and Alain Genestier on
the subject of these notes were very important. Particularly, they indicated to the second author
the analogy between Anderson varieties and unitary Shimura varieties,
and Alain Genestier attracted his attention to the paper [F] where the
notion of the dual Anderson module is defined, and indicated the statement of Theorem 3.4.1a. Finally, the second author is grateful to Nicholas Katz for a detailed explanation of his paper [K].
\endthanks
\address
First author: Departamento de Matem\'atica e estatistica
Universidade de S\~ao Paulo. Rua de Mat\~ao 1010, CEP 05508-090, S\~ao Paulo, Brasil, and Omsk State University n.a. F.M.Dostoevskii. Pr. Mira 55-A, Omsk 644077, Russia.
\medskip
Second author: Departamento de Matem\'atica, Universidade Federal do Amazonas, Manaus, Brasil
\endaddress
\abstract We formulate some properties of a conjectural object $X_{fun}(r,n)$ parametrizing Anderson t-motives of dimension $n$ and rank $r$. Namely, we give formulas for $\g p$-Hecke correspondences of $X_{fun}(r,n)$ and its reductions at $\g p$ (where $\g p$ is a prime of $\n F_q[\th]$). Also, we describe their geometric interpretation. These results are analogs of the corresponding results of reductions of Shimura varieties. Finally, we give conjectural formulas for Hodge numbers (over the fields generated by Hecke correspondences) of middle cohomology submotives of $X_{fun}(r,n)$.
\endabstract
\keywords Moduli objects of Anderson t-motives, Hecke correspondences, reduction \endkeywords
\subjclass 11G09 \endsubjclass
\endtopmatter

\document
{\bf 0. Introduction.} Let $X$ be a Shimura variety of PEL-type. Its points parametrize abelian varieties with some  fixed endomorphism rings, polarization and level. There is a problem to describe relations between the rings $\n H_p(X)$ of $p$-Hecke correspondences of $X$ ($p$ is a prime) and of $\tilde X_p$ - the reduction of $X$ at $p$, particularly, to find the characteristic polynomial of the Frobenius correspondence of $\tilde X_p$ over $\n H_p(X)$. Also, we can ask what are geometric interpretations of reductions of $p$-Hecke correspondences.
\medskip
There is also a problem of extreme complexity --- to prove Langlands theorems for Shimura varieties (relations between $L$-functions of submotives of $X$ and of automorphic representations of the corresponding reductive group $G$). It is solved only in a few cases of $X$ of low dimension. Even exact statements of theorems giving these relations is a complicated problem.
\medskip
Analogs of abelian varieties for the function field case are Anderson t-motives (or Anderson modules --- their categories are anti-equivalent). It is natural to consider an object (function field analog of a Shimura variety) whose points parametrize these t-motives (for example, Anderson t-motives of a fixed dimension $n$, rank $r$, type of a nilpotent operator $N$, see Section 2.4, and of an analog of PEL-type). We denote this object by $X_{fun}=X_{fun}(r,n)$.
\medskip
Unfortunately, at the moment for $n>1$ these objects are conjectural (only for the case $n=1$ we have a good theory of moduli spaces of Drinfeld modules). For example, if we restrict ourselves by pure uniformizable Anderson t-motives and assume that there is 1 -- 1 (or near 1 -- 1) correspondence between these t-motives and their lattices (which is rather likely, see [GL17]), then the moduli variety of lattices would be the quotient set of Siegel matrices by an (almost) action of $GL_r(\n F_q[\th])$. But this (almost) action does not have desired properties, see [GL17], Proposition 1.7.1.
\medskip
So, the whole contents of the present paper concerns conjectural objects $X_{fun}(r,n)$ that we shall call Anderson modular objects. Nevertheless, we can get some information about them. Let us give more definitions. Naively, $X_{fun}(r,n)$ parametrize pure abelian t-motives of rank $r$ and dimension $n$ whose nilpotent operator $N$ is 0. The corresponding reductive group $G$ is $GL_r$ and the dominant coweight $\mu$ is $(1,\dots, 1, 0,\dots, 0)$ ($n$ ones and $r-n$ zeroes). For more information for any $G$ on $X_{fun}$ and $\mu$ see [V].
\medskip
Let $\goth p$ be a prime ideal of $\n F_q[\th]$ ($q$ is a power of $p$ and $\th$ an abstract variable, see Section 2). The contents of the present paper is the following.
\medskip
(1) We formulate (for some cases) in Section 2 the theorems concerning reductions at $\goth p$ of Hecke
correspondences $T_{\goth p,i}$ ($i=0,\dots,r$) on $X_{fun}$, and their
geometric interpretation. These results are functional analogs of [FCh],
Chapter 7, [BR], Chapter 6, and [W]. The case $n=1$ (Drinfeld varieties) is treated with more details and
explicit formulas.
\medskip
(2) Sketches of the proofs of these theorems are given in Section 3.
\medskip
We see that the function field case --- Anderson varieties of rank $r$ and
dimension $n$, $G_{fun}=GL_r(\n F_q(\th))$ is analogous to the number field case where
$G_{num}=GU(r-n,n)$
corresponds to Shimura varieties (called unitary for brevity) of PEL-type parametrizing
abelian $r$-folds with multiplication by an imaginary quadratic field
$K$, of signature $(r-n,n)$. We indicate in Section 4 that really, properties of unitary Shimura varieties are similar to the properties of $X_{fun}(r,n)$. By the way, this analogy is a source of more results, see for example [GL21].
\medskip
(3) Finally, in Section 5 we state conjectural formulas for Hodge numbers (over fields generated by Hecke correspondences) of submotives of middle cohomology of $X_{fun}$. They are analogs of the corresponding formulas for Shimura varieties ([BR], SEction 4.3, p. 548). In Section 6 we consider the action of Hecke correspondences on
some non-ordinary Drinfeld modules. These results will be useful for a future
proof of analog of Kolyvagin's theorem (finiteness of Tate-Shafarevich group) for the
case of Drinfeld varieties.
\medskip
In order to show the analogy between the number field and the function fiels cases, we give in Section 1 some
well-known results on Hecke correspondences and their reductions for the case of Siegel varieties.
\medskip
{\bf 1. Definitions and results for the number field case (Siegel varieties).} \nopagebreak
\medskip
{\bf 1.1. Reductions of correspondences.} \nopagebreak
\nopagebreak
\medskip
For a comparison with the number field case, here we formulate well-known results for reductions of Siegel modular varieties. References for the results of this section are: [FCh], Section 7, and [W].
\medskip
Let $X$ be a Siegel variety of
genus $g$
(of any fixed level), i.e. a quotient of the Siegel upper half plane by a congruence subgroup of $GSp_{2g}(\n Z)$, or, equivalently, a set of principally polarized abelian varieties of dimension $g$ together with some level structure. Let the congruence subgroup be such that $X$ is defined over $\n Q$. We have $G=G_X=GSp_{2g}(\n Q)$ is the
corresponding reductive group. Let $p$ be a fixed prime which does not divide the
level, and $\tilde X$ the reduction of $X$ at $p$.
\medskip
Let $\Corr(X)=\Corr_p(X)$ (resp. $\Corr(\tilde X)=\Corr_p(\tilde X)$) be the
algebra of $p$-Hecke correspondences on $X$ (resp. $\tilde X$). We have the
Frobenius map on $\tilde X$; considering it as a correspondence we get an
element $fr_X \in \Corr(\tilde X)$.
\medskip
There is a map $\gamma : \Corr(X)  \to \Corr (\tilde X)$ --- the
reduction of a
correspondence at $p$. There is a problem of description of $\gamma$ and of
finding of the characteristic polynomial of $fr$ over $im(\gamma)$. The answer
is the following. Let $M$ be the following block diagonal subgroup of $G$:
$$M=\left\{\left( \matrix A & 0 \\ 0 & \lambda \cdot (A^t)^{-1} \endmatrix
\right)\right\}\subset G,\eqno{(1.1.1)}$$ (blocks have size $g$), and let $T$ be the
subgroup of $M$ consisting of diagonal matrices.
\medskip
The abstract $p$-Hecke algebras $\n H(G)=\n
H(G)(\n Q_p)$
(resp. $\n H(M)=\n H(M)(\n Q_p)$, $\n H(T)=\n H(T)(\n Q_p)$) consist of double
cosets $K\alpha K$, where $K=K_G=G(\n Z_p)$ (resp. $K=K_M=M(\n Z_p)$, $K=K_T=T(\n Z_p)$) and $\al\in G(\n Q_p)$, resp. $\al\in M(\n Q_p)$, $\al\in T(\n Q_p)$.
There are the Satake inclusions $S^G_M: \n H(G) \to \n H(M)$, $S^M_T: \n H(M) \to \n
H(T)$.
\medskip
The Hecke algebra $\n H(T)$ is the
subalgebra of $\n Z[U_1^{\pm 1},...,U_g^{\pm 1},V_1^{\pm 1},...,V_g^{\pm 1}]$ (here $U_i$, $V_i$ are abstract variables)
generated by $(U_1V_1^{-1})^{\pm 1},...,(U_gV_g^{-1})^{\pm 1}$ and
$(U_1\cdot...\cdot U_g)^{\pm 1}$. The Weyl group $W_G$ of $G$ is the
semidirect
product of the permutation group $S_g$ and of $(\n Z/2\n Z)^g=(\pm1)^g$, where $S_g$ permutes
coordinates in $(\pm1)^g$. There is a section $S_g\hookrightarrow W_G$, we
denote its image by $W_{G,M}$. Further, $W_G$ acts on $\n H(T)$ in the obvious
manner: $S_g$ permutes indices and $(\pm1)^g$ interchanges $U_*$, $V_*$.
\medskip
We have: $$S^M_T(\n H(M))=\n H(T)^{W_{G,M}}, \ \  S^G_T(\n H(G))=\n
H(T)^{W_G}.\eqno{(1.1.2)}$$
\medskip
It is known that there are surjections
$$\beta_1: \n H(G) \to \Corr(X),\eqno{(1.1.3)}$$ $$\beta_2: \n H(M) \to
\Corr(\tilde X)\eqno{(1.1.4)}$$
whose kernel is generated by the relation $K p K=id$, where $p=pI_{2g}$ is the scalar
matrix in both $G$ and $M$, and $id$ is the trivial correspondence on both $X$, $\tilde X$.
\medskip
{\bf Theorem 1.1.5.} There exists a commutative diagram:
$$\matrix S^G_M: & \n H(G) & \to & \n H(M) \\ \\
                       & \beta_1\downarrow & & \beta_2\downarrow \\ \\
                                  \gamma : &\Corr(X)  & \to & \Corr (\tilde X)
                                  \endmatrix \eqno{(1.1.6)}$$

We denote by $\tau_p$ the matrix $\left( \matrix 1 & 0 \\ 0 & p \endmatrix
\right) $, where entries are scalar $g \times g$-blocks, and we denote the
corresponding elements $K_G \tau_p K_G$ (resp. $K_M \tau_p K_M$) of $\n H(G)$
(resp. $\n H(M)$) by $T_p$ (resp. $fr_M$).
\medskip
{\bf Theorem 1.1.7.} $\beta_2(fr_M)=fr_X$.
\medskip
{\bf Remark.} Formulas 1.1.2 and theorems 1.1.5, 1.1.7 permit us to find the
Hecke polynomial of $X$ (the characteristic polynomial of $fr_X$ over $\Corr(X)$ ). Really, 1.1.2 implies that $\n H(M)$ is a free module
over $S^G_M(\n H(G))$ of dimension $\#(W_G)/\#(W_{G,M})=2^g$. An explicit
description of $fr_M \in \n H(M)$ (see below) permits us to find easily its
characteristic polynomial over $\n H(G)$. Theorem 1.1.5 implies that
it is also
the characteristic polynomial of the Frobenius correspondence on
$\tilde X$ over
the algebra $\Corr X$.
\medskip
{\bf 1.2. Geometric interpretation.}\nopagebreak
\medskip
For $ i = 0,\dots,g$ we consider diagonal matrices in a block form $$\varphi_i =
\left( \matrix
I_i & 0 & 0 & 0
                              \\ 0 & pI_{g-i} & 0 & 0
                              \\ 0 & 0 & pI_i & 0
                              \\ 0 & 0 & 0 & I_{g-i} \endmatrix
                                                          \right)\in M $$
(sizes of diagonal blocks are $i$, $g - i$, $i$, $g - i$), and we
denote the
corresponding elements $K_M \varphi_i K_M\in \n H(M)$ by $\Phi_i$.
Particularly,
$\Phi_g=fr_M$.
\medskip
Let $I$ be a subset of $\{1, ... ,g\}$. We denote $U_I:=\prod_{i\in I}U_i \prod_{i\not\in
I}V_i\in \n H(T)$. We have $$S^M_T(\Phi_i)=\sum_{\#(I)=i}U_I$$ and
$$S^G_M(T_p)=\Phi_0+\Phi_1+\dots +\Phi_g\in \n H(M).\eqno{(1.2.1)}$$
We denote $\beta_1(T_p) \in \Corr(X)$, $\beta_2(\Phi_i) \in \Corr(\tilde X)$
again by $T_p$, $\Phi_i$ respectively, so (1.2.1) and Theorem 1.1.5
give us the
following equality on $\Corr(\tilde X)$:
$$\gamma(T_p)=\Phi_0+\Phi_1+\dots +\Phi_g.\eqno{(1.2.2)}$$

For any algebraic variety $Z$ there exists an involution on $\Corr (Z)$ (symmetry with respect the coordinates). Also, there exist involutions on $ \n H(M)$, $ \n H(G)$ commuting with involutions on $\Corr (X)$, $\Corr (\tilde X)$ with respect to (1.1.6). We denote these involutions by hat; we have $\hat T_p=T_p$,  $\hat \Phi_i=\Phi_{g-i}$.
\medskip
The geometric interpretation of (1.2.2) is the following. Let
$t \in X(\bar \n Q)$ be a generic point, $A_t$ the corresponding principally polarized abelian
$g$-fold with a fixed polarization form and $(A_t)_p$ the $\n F_p$-space of its $p$-torsion points.
The polarization on $A_t$ defines a
skew form on $(A_t)_p$. $T_p(t)$ is a finite set of points; we have: $t'\in
T_p(t)$ iff
there exists an isogeny $\alpha_{t,t'}: A_t\to A_{t'}$ of type
$(1,...,1,p,...,p)$. The kernel of $\alpha_{t,t'}$ is an isotropic
$g$-dimensional subspace of $(A_t)_p$. So, we have a
\medskip
{\bf Theorem 1.2.3.} The set $T_p(t)$ is in 1--1 correspondence with
the set of
isotropic $g$-dimensional subspaces of $(A_t)_p$.
\medskip
Now let $t\in X(\n Q)$ be a generic point such that $A_t$ has a good ordinary reduction at $p$. Let $(\tilde A_t)_{p,points}$ be the set of closed points of $\tilde
A_t$ of order $p$, and $red: (A_t)_p \to (\tilde A_t)_{p,points}$ the
reduction map. We denote by $D=D_{Siegel}$ the kernel of $red$, it is an
isotropic $g$-dimensional subspace of $(A_t)_p$.
\medskip
Let $t'\in T_p(t)$ and $\tilde t'\in \tilde X$ its reduction. (1.2.2)
shows that
$\tilde t'$ belongs to one of $\Phi_i(\tilde t)$.
\medskip
{\bf Theorem 1.2.4.} Number $i$ is defined as follows:
$$i=\dim(\Ker(\alpha_{t,t'})\cap D_{Siegel}).\eqno{(1.2.5)}$$
Particularly, $i=g \iff \tilde t'=fr(\tilde t)\iff \Ker(\alpha_{t,t'})=
D_{Siegel}$.
\medskip
Further, we have the following
\medskip
{\bf Theorem 1.2.6.} Let $t'$, $t''$ be 2 points of $T_p(t)$. Then
$$\tilde t'=\tilde t''\iff \Ker(\alpha_{t,t'})\cap D_{Siegel} =
\Ker(\alpha_{t,t''})\cap D_{Siegel}.$$

Recall that any correspondence $C$ on $X$ has the bidergee $d_1(C),d_2(C)$ ---
the degrees of 2 projections $\pi_1$, $\pi_2$ of its graph $\Gamma(C)$ to $X$. By definition, $\pi_{i}(\hat C)=\pi_{\hat i}(C)$ (here $\hat 1=2, \ \hat 2=1$).
Further, $C$
has the separable (resp. non-separable) bidergee
$d^{s}_1(C),d^{s}_2(C)$ (resp.
$d^{ns}_1(C),d^{ns}_2(C)$) --- the separable (resp. non-separable) degrees of
$\pi_1$, $\pi_2$. We have $d_i(C)=d^{s}_i(C)d^{ns}_i(C)$ and $$d_{i}^*(\hat C)=d_{\hat i}^*(C), \ \ \ i=1,2, \ \ \ \ *=\emptyset, s, ns.\eqno{(1.2.6a)}$$ We
denote by
$\goth g(k,l)$ the cardinality of the Grassmann variety $Gr(k,l)(\n F_p)$:
$$\goth g(k,l)=\prod_{i=1}^k \frac{p^l-p^{i-1}}{p^k-p^{i-1}}.\eqno{(1.2.7)}$$
{\bf Theorem 1.2.8.}
$$d_1^{s}(\Phi_i)=\goth g(i,g), \ d_1^{ns}(\Phi_i)=p^{(g+1-i)(g-i)/2}, \
d_2^{s}(\Phi_i)=\goth g(i,g), \ d_2^{ns}(\Phi_i)=p^{(i+1)i/2}. \ \square$$

{\bf 2. Definitions and statement of conjectures for the case of
Anderson modular objects.}\nopagebreak
\medskip
We use standard notations for Anderson t-motives. Let $q$ be a power of a prime $p$, $\n F_q$ the finite field of order $q$. The function field analog of $\n Z$ is the ring of polynomials $\n F_q[\th]$ where $\th$ is an abstract variable. The analog of the archimedean valuation on $\n Q$ is the valuation at infinity of the fraction field $\n F_q(\th)$ of $\n F_q[\th]$; it is denoted by $ord$, it is uniquely determined by the property $\ord(\th)=-1$. The completion of an algebraic closure of the completion of $\n F_q(\th)$ with respect the valuation "ord" is the function field analog of $\n C$. It is denoted by $\p$.
\medskip
The definition of a t-motive $\g M$ is given in [G], Definitions 5.4.2,
5.4.12 (Goss uses another terminology: "abelian t-motive" of [G] = "t-motive" of the present paper). Particularly, $\g M$ is a free $\p[T]$-module of
dimension $r$ (this number $r$ is called the rank of $\g M$) endowed by a $\p$-skew-linear
operator $\tau$ satisfying some properties. Its dimension $n$
is defined in [G], Remark 5.4.13.2 (Goss denotes the dimension by $\rho$). A nilpotent operator $N=N(\g M)$
associated to a t-motive is defined in [G], Remark 5.4.3.2. Condition
$N=0$ implies $n \le r$. Except Section 2.4, we shall consider only the case $N=0$.
\medskip
As it was written in the introduction, the main object of the present paper is conjectural. It is called an Anderson modular object, it is denoted by $X_{fun}=X_{fun}(r,n)$, it is the function field analog of $X$. Naively, it parametrizes Anderson t-motives of rank $r$ and dimension $n$.
\medskip
An analog of $p$ of Section 1 is a valuation (distinct of ord) of $\n F_q(\th)$ = a prime ideal of $\n F_q[\th]$. We denote by $\goth p$ both its generator and
the prime ideal itself, and we denote $\goth q = \#(\n F_q[\th]/\goth p)$. The
corresponding algebraic group $G_{fun}$ --- the function field analog of
$GSp_{2g}(\n Q)$ --- is $GL_r(\n F_q(\th))$. Hence, the analogs of $\n Q_p$, $\n Z_p$ for the functional case are $\n F_q(\th)_\g p$, $\n F_q[\th]_\g p$ respectively, and the analog of $K_G$ of Section 1 is $GL_r(\n F_q[\th]_\g p)$ (it will be denoted by $K_G$ as well).
\medskip
In order to simplify the present version of the text, for the case $n>1$ we
state conjectures of Section 2.3 only for uniformizable Anderson
t-motives. Analytically, an uniformizable Anderson t-motives of rank $r$ and dimension $n$
over $\p$ is the quotient $\p/L$, where $L$ is a free
$r$-dimensional $\n F_q[\th]$-module.  Since not all Anderson t-motives are
uniformizable, the exact statements of these conjectures must be slightly
changed, see Remark 2.3.4a.

\medskip
{\bf Theorem 2.1.} The analog of $M$ for this case is the group $M_{r-n,n}$
of block diagonal matrices $\left( \matrix *&0\\0&* \endmatrix \right)\subset
GL_r(\n F_q(\th))$, sizes of blocks are $r-n$, $n$. We denote $K_M:=M(\n F_q[\th]_\g p)$. The analogs of Hecke algebras and of algebras of correspondences are defined like in Section 1. The analogs of Theorems 1.1.5, 1.1.7
hold for this case. Particularly, $\Corr(\tilde X_{fun})$ is the quotient of
$\n H(M)$ by the trivial relation $K_M\goth p K_M=id$.
\medskip
{\bf 2.2. Description of $\n H(G_{fun})$, $\n H(M)$ and of the Satake
inclusions.}\nopagebreak
\medskip
Let like in Section 1, $T$ be the subgroup of $G_{fun}$ of diagonal matrices. We have
$\n H(T)=\n Z[U_1^{\pm 1},...,U_r^{\pm 1}]$, and the Weyl group
$W_{G_{fun}}=S_r$, it acts on $\n H(T)$ permuting indices. An analog of
$W_{G,M}$ is the subgroup $W_{G_{fun},M}=S_{r-n}\times S_n\hookrightarrow
S_r=W_{G_{fun}}$ with the obvious inclusion. Formulas (1.1.2) hold for our case,
explicit formulas are the following.
\medskip
{\bf 2.2.1.} $\n H(G_{fun})$: For $ i = 0,\dots,r$ we denote by $\tau_{\goth
p,i}$ the diagonal matrix $\left( \matrix I_{r-i} & 0 \\ 0 & \goth pI_i \endmatrix
\right) \in G_{fun}$, where sizes of blocks are $r - i$, $i$, and we
denote the
corresponding elements $K_G \tau_{\goth p,i} K_G\in \n H(G_{fun})$ by
$T_{\goth
p,i}$. We have $T_{\goth p,0}=1$, $T_{\goth p,r}$ is the trivial correspondence, and other $T_{\goth p,i}$ are free
generators
of $\n H(G_{fun})$. We have $$S^{G_{fun}}_T(T_{\goth p,i})=\goth
q^{-i(i-1)/2}\sigma_i(U_1,...,U_r),\eqno{(2.2.2)}$$ where $\sigma_i$ is the
$i$-th symmetric polynomial, and $$\hat T_{\goth p,i}=T_{\goth p,r-i}.\eqno{(2.2.2a)}$$
\medskip
{\bf 2.2.3.} $\n H(M)$: (a) For $ i = 0,\dots,r-n$ we denote by
$\varphi_i$ the
diagonal matrix
$\left( \matrix I_{r-n-i} & 0 & 0 \\ 0 & \goth pI_i & 0 \\ 0 & 0 & I_n \endmatrix
\right)\in
M $
where sizes of blocks are $r-n - i$, $i$, $n$, and we denote the corresponding
elements $K_M \varphi_i K_M\in \n H(M)$ by $\Phi_i$.

(b) For $ i = 0,\dots,n$ we denote by $\psi_i$ the diagonal matrix
$\left( \matrix I_{r-n} & 0 & 0 \\ 0 & \goth pI_i & 0 \\ 0 & 0 & I_{n-i} \endmatrix
\right)\in
M $
where sizes of blocks are $r-n$, $i$, $n-i$, and we denote the corresponding
elements $K_M \psi_i K_M\in \n H(M)$ by $\Psi_i$. We have $\Phi_0=\Psi_0=1$,
$\Phi_i=0$ (resp. $\Psi_i=0$) if $i$ is out of the range $0,\dots,r-n$ (resp.
$0,\dots,n$) and other $\Phi_i$, $\Psi_i$ are free generators of $\n H(M)$.
Obviously
$$S^{M}_T(\Phi_i)=\goth q^{-i(i-1)/2}\sigma_i(U_1,...,U_{r-n}), \
S^{M}_T(\Psi_i)=\goth
q^{-i(i-1)/2}\sigma_i(U_{r-n+1},...,U_{r}).\eqno{(2.2.4)}$$
\medskip
Formulas (2.2.2), (2.2.4) imply immediately that
$$S^{G_{fun}}_M(T_{\goth p,j})=\sum_{i=0}^j \goth
q^{-i(j-i)}\Psi_i\Phi_{j-i}.\eqno{(2.2.5)}$$
Further, we have:
$$\Psi_n\Phi_{r-n}=\goth q^{n(r-n)},\eqno{(2.2.5a)}$$
$$\hat \Phi_i=\goth q^{-n(r-n-i)}\Psi_n\Phi_{r-n-i},\eqno{(2.2.5b)}$$
$$\hat \Psi_i=\goth q^{-(n-i)(r-n)}\Psi_{n-i}\Phi_{r-n}.\eqno{(2.2.5c)}$$
particularly, $\hat \Phi_{r-n}=\Psi_n$. Coefficients of (2.2.5b,c) can be found from the property that equations (2.2.2a) and (2.2.5) are concordant with respect to the duality.
\medskip
{\bf Remark.} (2.2.5), (2.1), (1.1.5) imply that for the case $n=1$ (Drinfeld modules) the explicit formulas
are the following ($fr=\Psi_1$):
$$\tilde T_{\goth p,1}=fr+\Phi_1;$$
$$\tilde T_{\goth p,2}=\frac1{\goth q} fr\ \Phi_1+\Phi_2;$$
$$\dots\eqno{(2.2.6)}$$
$$\tilde T_{\goth p,r-1}=\frac1{\goth q^{r-2}}fr\ \Phi_{r-2}+\Phi_{r-1};$$
$$\tilde T_{\goth p,r}=\frac1{\goth q^{r-1}}fr\ \Phi_{r-1}.$$

$\n H(M)$ is a free module over $S^{G_{fun}}_M(\n H(G_{fun}))$
respectively the
Satake inclusion. Its dimension is $\#(W_{G_{fun}})/\#(W_{G_{fun},M})=\left(
\matrix r \\ n \endmatrix \right)$.

We denote $\Psi_n$ by $fr_M$. Its characteristic polynomial over
$S^{G_{fun}}_M(\n H(G_{fun}))$ (the Hecke polynomial) can be easily found by
elimination of $\Phi_1,...,\Phi_{r-n},\Psi_1,...,\Psi_{n-1}$ in the system
(2.2.5). We denote it by $P_{r,n}$, it belongs to $\n H_{\goth
p}(G_{fun})[fr_M]$. For $n=1$ we have
$$\matrix P_{r,1}=\sum_{i=0}^r (-1)^i \goth q^{i(i-1)/2} T_{\goth p,i}
fr^{r-i}
\\  \\ =fr^r - T_{\goth p,1}fr^{r-1}+\goth q T_{\goth p,2}fr^{r-2}\pm\dots
+(-1)^r\goth q^{r(r-1)/2}T_{\goth p,r}.\endmatrix \eqno{(2.2.7)}$$

{\bf 2.3. Statements of results. }\nopagebreak
\medskip
Theorems 1.2.4, 1.2.6 can be rewritten almost word-to-word as conjectures for
the functional case. Let us do it. It is more convenient to use Anderson modules ([G], 5.4.5; Goss calls them Anderson $T$-modules) instead of Anderson t-motives. The categories of Anderson t-motives and modules are anti-isomorphic, so there is no essential difference which object to use.

Let $t\in X_{fun}$ be such that the
corresponding Anderson module $E_t$ is uniformizable: $E_t=\p^n/L$ (as earlier, $n$ is the dimension and $r$ is the rank). This condition is ``closed under Hecke correspondences'': if $t'\in
T_{\goth p,j}(t)$ then $E_{t'}$ is also uniformizable. Obviously, we have
\medskip
{\bf Theorem 2.3.1.} $(E_t)_{\goth p}$ --- the set of ${\goth
p}$-torsion points
of $E_t$ --- is $\goth p^{-1}L/L$ and hence is an $r$-dimensional $\n F_q[\th]/\goth
p$-vector space.
\medskip
{\bf Theorem 2.3.2.} $t'\in T_{\goth p,j}(t)$ iff there exists an isogeny
$\alpha_{t,t'}: E_t\to E_{t'}$ of type $(1,...,1,\goth p,...,\goth p)$ ($r-j$
1's and $j$ ${\goth p}$'s).
\medskip
{\bf Theorem 2.3.3.} The set $T_{\goth p,j}(t)$ is in 1--1 correspondence with
the set of $j$-dimensional subspaces of $(E_t)_{\goth p}$.
\medskip
Now we can formulate the conjecture on reductions at ${\goth p}$. Let $(\tilde
E_t)_{\goth p,points}$ and $red: (E_t)_{\goth p} \to (\tilde E_t)_{\goth
p,points}$ be analogs of the corresponding objects in Section 1.2. We
denote by
$D_{fun}$ the kernel of $red$. Firstly, we have a
\medskip
{\bf Theorem 2.3.4.} For a generic $t\in X_{fun}$ \ \ $D_{fun}$ is an
$n$-dimensional subspace of $(E_t)_{\goth p}$.
\medskip
{\bf Remark 2.3.4a.} For an arbitrary $t\in X_{fun}$ (such that $E_t$ is not
uniformizable) statements of the above and below theorems and conjectures
require minor changes. For example, Theorem 2.3.1 becomes
\medskip
{\bf Theorem 2.3.4b.} $(E_t)_{\goth p}$ is an $r$-dimensional $A/\goth
p$-vector space.
\medskip
(we cannot claim that it is $\goth p^{-1}L/L$ because $L$ does not exist).
\medskip
{\bf Geometric interpretation.}

\nopagebreak
\medskip
Here we formulate analogs of (1.2.4) -- (1.2.8) for the function field case. Let
$t'\in T_{\goth p,j}(t)$ and $\tilde t'\in \tilde X_{fun}$ its reduction.
Conjecture 2.1 and (2.2.5) show that $\tilde t'$ belongs to one of $\goth
q^{-i(j-i)}\Psi_i\Phi_{j-i}(\tilde t)$.
\medskip
{\bf Theorem 2.3.5.} Number $i$ is defined as follows:
$$i=\dim(\Ker(\alpha_{t,t'})\cap D_{fun}).\eqno{(2.3.6)}$$

Particularly, $j=i=n \iff \tilde t'=fr_{X_{fun}}(\tilde t)\iff
\Ker(\alpha_{t,t'})= D_{fun}$.
\medskip
Now $\goth g(k,l)$ will mean the cardinality of Grassmann variety
$Gr(k,l)$ over
$\n F_q[\th]/\goth p$, i.e. $p$ in (1.2.7) must be replaced by $\goth q$. Obviously we
have
$$d_1(T_{\goth p,i})=d_2(T_{\goth p,i})=\goth g(i,r).$$

For the reader's convenience, we formulate the following conjecture separately
for the case $n=1$.
\medskip
{\bf Theorem 2.3.7.} For $n=1$ we have:
$$d_1^{s}(\Psi_1)=d_1^{ns}(\Psi_1)=1, \ d_2^{s}(\Psi_1)=1, \
d_2^{ns}(\Psi_1)=\goth q^{r-1};$$
$$d_1^{s}(\Phi_i)=\goth g(i,r-1), \ d_1^{ns}(\Phi_i)=\goth q^i, \
d_2^{s}(\Phi_i)=\goth g(i,r-1), \ d_2^{ns}(\Phi_i)=1.$$

{\bf Corollary 2.3.7.1.} For correspondences $\frac1{\goth q^{i}}fr\ \Phi_{i}$
--- summands in the right hand side of (2.2.6) --- we have
$$d_1^{s}(\frac1{\goth q^{i}}fr\ \Phi_{i})=\goth g(i,r-1), \ \
d_1^{ns}(\frac1{\goth q^{i}}fr\ \Phi_{i})=1, $$ $$ d_2^{s}(\frac1{\goth
q^{i}}fr\ \Phi_{i})=\goth g(i,r-1), \ \ d_2^{ns}(\frac1{\goth q^{i}}fr\
\Phi_{i})=\goth q^{r-1-i}.$$

An analog of the Theorem 1.2.6 is the following. Let $t'$, $t''$ be 2 points of
$T_{\goth p,i}(t)$. Firstly we consider the case when
$$\Ker(\alpha_{t,t'})\cap D_{fun}=\Ker(\alpha_{t,t''})\cap D_{fun}=0,$$

i.e. both $\tilde t'$, $\tilde t''\in \Phi_i(\tilde t)$ (the second summand in
(2.2.6)).
\medskip
{\bf Theorem 2.3.8.} Under this condition we have (for any $n$):
\medskip
$\tilde t'=\tilde t''$ as closed points iff the linear spans coincide:
$$<\Ker(\alpha_{t,t'}), D_{fun}>=<\Ker(\alpha_{t,t''}), D_{fun}>.$$
If (for $n=1$)
$$\Ker(\alpha_{t,t'})\supset D_{fun}, \ \Ker(\alpha_{t,t''})\supset D_{fun},$$
i.e. both $\tilde t'$, $\tilde t''\in \frac1{\goth
q^{i-1}}\Psi_1\Phi_{i-1}(\tilde t)=\frac1{\goth q^{i-1}}fr \ \Phi_{i-1}(\tilde
t)$ (the first summand in (2.2.6)) then (2.3.7.1) implies that $\tilde t'$,
$\tilde t''$ are always different.
\medskip
Now let us consider the case of arbitrary $n$.
\medskip
{\bf Theorem 2.3.10.}
$$d_1^{s}(\Psi_i)=\goth g(i,n), \ d_1^{ns}(\Psi_i)=1;$$
$$d_1^{s}(\Phi_i)=\goth g(i,r-n), \ d_1^{ns}(\Phi_i)=q^{in}.$$

{\bf Remark 2.3.10a.} Numbers $d_2^*(\Phi_i)$, $d_2^*(\Psi_i)$ ($*=\emptyset, s, ns$) can be found from the above formulas using (1.2.6a), (2.2.5b,c).  Particularly, for the Frobenius $\Psi_n$ we have
$$d_2^{s}(\Psi_n)=1, \ d_2^{ns}(\Psi_n)=q^{n(r-n)}.$$
\medskip
{\bf Corollary 2.3.11.} For correspondences $\frac1{\goth
q^{i(j-i)}}\Psi_i\Phi_{j-i}$ --- summands in the right hand side of
(2.2.5) ---
we have
$$d_1^{s}(\frac1{\goth q^{i(j-i)}}\Psi_i\Phi_{j-i})=\goth g(i,n)\goth
g(j-i,r-n), \ d_1^{ns}(\frac1{\goth q^{i(j-i)}}\Psi_i\Phi_{j-i})=\goth
q^{(j-i)(n-i)}.$$

Let $t'$, $t''\in T_{\goth p,j}(t)$ be as above such that
$$\dim(\Ker(\alpha_{t,t'})\cap D_{fun})=\dim(\Ker(\alpha_{t,t''})\cap
D_{fun})=i.$$
According Theorem 2.3.5, both $\tilde t'$, $\tilde t''\in \frac1{\goth
q^{i(j-i)}}\Psi_i\Phi_{j-i}$.
\medskip
{\bf Theorem 2.3.12.} $\tilde t'=\tilde t''$ as closed points iff both
intersections and linear spans coincide:
$$\Ker(\alpha_{t,t'})\cap D_{fun}= \Ker(\alpha_{t,t''})\cap D_{fun};$$
$$<\Ker(\alpha_{t,t'}), D_{fun}>=<\Ker(\alpha_{t,t''}), D_{fun}>.$$

{\bf Remark 2.3.13.} (a) If we fix an $i$-dimensional subspace $V_i \subset
D_{fun}$ and a $(n+j-i)$-dimensional overspace $V_{n+j-i} \supset D_{fun}$ then
the quantity of $j$-dimensional spaces $\alpha$ such that $$\alpha\cap
D_{fun}=V_i, \ \ <\alpha,D_{fun}>=V_{n+j-i},$$ is equal to $\goth
q^{(j-i)(n-i)}=d_1^{ns}(\frac1{\goth q^{i(j-i)}}\Psi_i\Phi_{j-i})$, as it is
natural to expect.
\medskip
(b) Thanks to existence of skew pairing in the number case, we have
$$<\Ker(\alpha_{t,t'}), D_{fun}>=(\Ker(\alpha_{t,t'})\cap D_{fun})^{dual},$$
($dual$ is with respect to the skew pairing), so in Theorem 1.2.6 it is sufficient to claim only coincidence of
intersections.
\medskip
(c) For $n=1$ (Theorem 2.3.8) intersections always coincide, so it is
sufficient to claim only coincidence of linear spans.
\medskip
{\bf 2.4. Case of} $N\ne0$. This is a subject of further research. Here we do not even give statements of results, we indicate only the discrete invariants of t-motives having $N\ne0$ and explain the methods how these statements can be obtained.
\medskip
Let $\g M$ be a uniformizable Anderson t-motive such that its $N$ is not (necessarily) 0, of dimension $n$ and of rank $r$. $\g M$ is a $\p[T]$-module with a skew map $\tau: \g M\to \g M $ such that $\g M/\tau \g M$ is annihilated by a power of $T-\theta$. We have $\Lie(\g M)=\p^n$, and $N$ is a nilpotent operator acting on $\Lie(\g M)$. Also, $T$ acts on $\Lie(\g M)$; we have $T=\th I_n+N$ on $\End(\Lie(\g M))$. Particularly, we can consider $\Lie(\g M)$ as a $\n F_q[T]$-module.
\medskip
The lattice $L(\g M)$ of $\g M$ is a free $\n F_q[T]$-submodule of $\Lie(\g M)$ considered as a $\n F_q[T]$-module. We have a natural inclusion $\n F_q[T]\hookrightarrow \p[[N]]$ ($T\mapsto N+\th$). Hence, the tautological inclusion $L(\g M) \hookrightarrow \Lie(\g M)$ defines a surjection
$$L(\g M)\underset{\n F_q[T]}\to{\otimes}\p[[N]]\twoheadrightarrow \Lie(\g M).$$

Its kernel is denoted by $\g q_\g M$; the exact sequence
$$0\to\g q_\g M\to L(\g M)\underset{\n F_q[T]}\to{\otimes}\p[[N]]\to \Lie(\g M)\to0$$
(see [P]; [Gl18], (1.8.3); [HJ], Example 2.5, (2.2) for the case of Drinfeld modules) is a particular case of a Hodge-Pink structure.
\medskip
Discrete invariants of $\g M$ are the discrete invariants of the pair of lattices $(\g q_\g M; L(\g M)\underset{\n F_q[T]}\to{\otimes}\p[[N]])$ over a discrete valuation ring $\p[[N]]$. Let us give a description of these invariants from [GL18], (3.3). Let $\nu$ be the minimal number such that $N^\nu=0$. These invariants are numbers $k_1\ge0, \dots, k_{\nu+1}\ge0$ defined as follows. There exists a basis $l_1, \dots, l_r$ of $L(\g M)$ over $\n F_q[T]$ and its partition on $\nu+1$ sets of lengths $k_1, \dots, k_{\nu+1}$ (if some $k_i=0$ then the $i$-th set is empty); the $i$-th set is denoted by $l_{i,1}\dots,l_{i,k_i}$, having the following properties:
\medskip
$N^{\nu-1}(l_{\nu+1,i})$, $i=1, \dots, k_{\nu+1}$, form a $\p$-basis of $N^{\nu-1}\Lie(\g M)$ ([GL18], (3.6);
\medskip
\medskip
$N^{\nu-2}(l_{\nu,i})$, $i=1, \dots, k_{\nu}$, $N^{\nu-2}(l_{\nu+1,i})$, $i=1, \dots, k_{\nu+1}$,
\medskip
$N^{\nu-1}(l_{\nu+1,i})$, $i=1, \dots, k_{\nu+1}$, form a $\p$-basis of $N^{\nu-2}\Lie(\g M)$ ([GL18], (3.8);
\medskip
\medskip
$N^{\nu-3}(l_{\nu-1,i})$, $i=1, \dots, k_{\nu-1}$, $N^{\nu-3}(l_{\nu,i})$, $N^{\nu-2}(l_{\nu,i})$, $i=1, \dots, k_{\nu}$,
\medskip
and $N^{\nu-3}(l_{\nu+1,i})$, $N^{\nu-2}(l_{\nu+1,i})$, $N^{\nu-1}(l_{\nu+1,i})$, $i=1, \dots, k_{\nu+1}$,
\medskip
form a $\p$-basis of $N^{\nu-3}\Lie(\g M)$ ([GL18], (3.8);
\medskip
\medskip
etc., until
\medskip
$l_{2,i}, \ i=1, \dots, k_{2},\  \dots, \ N^{\nu-1}(l_{\nu+1,i})$, $i=1, \dots, k_{\nu+1}$,
\medskip
form a $\p$-basis of $\Lie(\g M)$.
\medskip
See [GL18], (3.3) - (3.10) for more details.
\medskip
Particularly, we have: $$r=k_1+...+k_{\nu+1}; \ \ n=k_2+2k_3+3k_4+...+\nu k_{\nu+1}.$$ ([GL18], (3.10) and (3.5)). If $\nu=1$, i.e. $N=0$ --- this is the case considered above, then $k_1=r-n, \ k_2=n$. Therefore, numbers $k_1, \dots,k_{\nu+1}$ are $N\ne0$-generalizations of numbers $r-n,\ n$.
\medskip
{\bf Conjecture 2.4.1.} The analog of the subgroup $M$ of $G_{fun}=GL_r$ for the set of Anderson t-motives having $N\ne0$ and invariants $k_1, \dots,k_{\nu+1}$ is the subgroup of $GL_r$ of block diagonal matrices with block sizes $k_1, \dots,k_{\nu+1}$.
\medskip
{\bf Remark 2.4.2.} Some of $k_i$ can be 0. In this case we cannot distinguish between $M$ for different sets of $k_1, \dots,k_{\nu+1}$. Hence, maybe it is necessary to modify the statement of Conjecture 2.4.1.
\medskip
As it was written, finding of analogs of statements of Sections 2.2, 2.3 for the sets of Anderson t-motives having invariants $k_1, \dots,k_{\nu+1}$ is a subject of further research.
\medskip
{\bf Remark 2.4.3.} Hartl and Juschka use some other invariants of $\g M$, see [HJ], Section 2. First, they consider slightly more general objects, namely, their $\g q=\g q_\g M$ is a subset not of $ L(\g M)\underset{\n F_q[T]}\to{\otimes}\p[[N]]$ but of $ L(\g M)\underset{\n F_q[T]}\to{\otimes}\p((N))$ (also, they consider a weight filtration on $L(\g M)$ ). Further, their Hodge-Pink weights $\om_1,\dots,\om_r$ are related with $k_1,\dots,k_{\nu+1}$ as follows: for all $i=1,\dots, \nu+1$ the number $-i+1$ occurs $k_i$ times among $\om_1,\dots,\om_r$ (i.e. among $\om_1,\dots,\om_r$ there are $k_1$ zeroes, $k_2$ minus ones etc.).
\medskip
\medskip
{\bf 3. Proofs.} We follow [FCh], Ch. 7, Section 4 using the same notations if
possible, and indicating results that are not completely analogous to the
number field case.
\medskip
Recall that $\g p$ is a prime ideal of $\n F_q[\th]$. We denote by $\n F_q[\th]_\g p$, $\n F_q(\th)_\g p$ the completions at $\g p$ of $\n F_q[\th]$, $\n F_q(\th)$ respectively, and by $\n F_q[\th]_\g p^{nr}$ the ring of integers of the maximal unramified extension of $\n F_q(\th)_\g p$. As usual, bar means an algebraic closure.
\medskip
There are maps $\n F_q[\th]_\g p^{nr}\hookrightarrow\overline{\n F_q(\th)_\g p}$, $\n F_q[\th]_\g p^{nr}\twoheadrightarrow\overline{\n F_q[\th]/\g p}$. The corresponding maps of schemes $\Spec \overline{\n F_q(\th)_\g p} \to \Spec \n F_q[\th]_\g p^{nr}$, $\Spec \overline{\n F_q[\th]/\g p} \to \Spec \n F_q[\th]_\g p^{nr}$ are denoted by $\xi_k$, $\xi_{\goth p}$
respectively. The inverse image $\xi_{\goth p}^*$ of an object (i.e. the
reduction of this object) is denoted by tilde.
\medskip
We fix $i$, and let $\Gamma$ be the graph of $T_{\goth p,i}$ over $\Spec
\n F_q[\th]_{\goth p}^{nr}$.
\medskip
It is known that it exists. For $t\in \Gamma$ (resp. $t\in \tilde \Gamma$) let
$\phi_t: E_t\to E'_t$ be the corresponding map of Anderson modules over $\Spec
\n F_q[\th]_{\goth p}^{nr}$ (resp. $\Spec \overline{\n F_q[\th]/{\goth p}}$).
\medskip
We consider the ordinary locus $\Gamma^0$ of $\Gamma$:

$$t\in \Gamma^0 \iff \xi_k(E_t), \ \xi_k(E'_t)$$
are ordinary.
\medskip
{\bf Lemma 3.1.} $\tilde \Gamma^0$ is dense in $\tilde \Gamma$. $\square$
\medskip
{\bf 3.1a.} Now let $\tau_{\goth p}\in G_{fun}$ be any diagonal matrix,
$T_{G,\goth p}$  the element of Hecke algebra $\n H(G_{fun})$ corresponding to
the double coset $K_G\tau_{\goth p}K_G$, $\Gamma$ the graph of $T_{G,\goth p}$
over $\Spec \n F_q[\th]_{\goth p}^{nr}$, and $\tilde \Gamma^0$ for this $\Gamma$ is
defined as earlier.
\medskip
Let $c$ be the highest power of ${\goth p}$ that appears in the
diagonal entries
of $\tau_{\goth p}$ (for example, if $\tau_{\goth p}=\tau_{\goth p,i}$ then
$c=1$). Let $s\in \tilde \Gamma^0$ and $E_s$, $E'_s$ the corresponding Anderson t-motives over $\Spec
\overline{\n F_q[\th]/{\goth p}}$. This means that we have a direct sum
decomposition of the finite $\n F_q[\th]$-module scheme $(E_s)[{\goth p}^c]$ over $\Spec
\overline{\n F_q[\th]/{\goth p}}$ on its multiplicative and etale part:

$$(E_s)[{\goth p}^c]=(E_s)[{\goth p}^c]_{mult}\oplus (E_s)[{\goth
p}^c]_{et},\eqno{(3.1.1)}$$
where
$$(E_s)[{\goth p}^c]_{mult}=(\mu_{\goth p^c})^n,\eqno{(3.1.2)}$$
$$(E_s)[{\goth p}^c]_{et}=(\Spec \n F_q[\th]/{\goth p}^c)^{r-n}.\eqno{(3.1.3)}$$

We can restrict $\phi_s$ to $(E_s)[{\goth p}^c]$ getting a map
$$(\phi_s)[{\goth p}^c]: (E_s)[{\goth p}^c] \to (E'_s)[{\goth
p}^c].\eqno{(3.1.4)}$$

In its turn, this map is restricted to both etale and multiplicative parts:
$$(\phi_s)[{\goth p}^c]_{mult}: (E_s)[{\goth p}^c]_{mult} \to (E'_s)[{\goth
p}^c]_{mult}\eqno{(3.1.5)}$$ and
$$(\phi_s)[{\goth p}^c]_{et}: (E_s)[{\goth p}^c]_{et} \to (E'_s)[{\goth
p}^c]_{et}.\eqno{(3.1.6)}$$
Taking into consideration (3.1.2), (resp, (3.1.3)), we see that
$\phi_s$ defines
elements in $\n H(GL_n)$ (resp. $\n H(GL_{r-n})$). In
concordance of notations of [FCh], we denote them by $a$ (resp. $d$).
This pair
($a$, $d$) defines us an element of $\n H(M)$. It is called the type of $s$.
\medskip
{\bf Remark.} Unlike in the number case, here the elements $a$, $d$ are independent.
\medskip
In order to formulate the below proposition 3.4, we need the following
notations:
\medskip
{\bf 3.2.} Let $\delta: \goth E_1 \to \goth E_2$ be a map of Anderson
modules over $\Spec \n F_q[\th]/{\goth p}$ of type $\Psi_i\Phi_j$.
\medskip
This means that $c$ of (3.1.1) is 1, and kernels of the map (3.1.5) (resp.
(3.1.6)) is isomorphic to $(\mu_{\goth p})^i$ (resp. $(\Spec \n F_q[\th]/{\goth
p})^{j}$). We denote them by $\tilde K_m$, $\tilde K_e$ respectively.
\medskip
Further, let $E_1$ be a Anderson module over $\Spec \n F_q[\th]_{\goth p}$ such
that $\tilde E_1=\goth E_1$.
\medskip
{\bf Lemma 3.3.} We can identify $\tilde K_m$ (resp. $\tilde K_e$)
with some $i$
(resp. $j$)-dimensional subspaces in $D_{fun}(E_1)$ (resp. $(E_1)_{\goth
p}/D_{fun}(E_1)$ ).
\medskip
We denote these subspaces by $K_m$, $K_e$ respectively.
\medskip
Now let us consider the set of pairs $(\phi, E_2)$ where $\phi: E_1
\to E_2$ is
a map of  Anderson modules over $\Spec \n F_q[\th]_{\goth p}$, such that $\tilde
\phi=\delta$ (and hence $\tilde E_2=\goth E_2$).
\medskip
{\bf Proposition 3.4.} The set of the above $(\phi, E_2)$ is isomorphic to the
set of subspaces $\Cal W\subset (E_1)_{\goth p}$ such that
$$\Cal W\cap D_{fun}(E_1)=K_m, \ \ \Cal
W+D_{fun}(E_1)/D_{fun}(E_1)=K_e.\eqno{(3.4.1)}$$
\medskip
{\bf Proof.} We need the function field analog of [K], Th. 2.1. Let $R$ be
an Artinian
local ring with residue field $\overline{\n F_q[\th]/{\goth p}}$ and the maximal ideal
$\goth m$. We consider only the case $R=R_\eta=\n F_q[\th]_{\goth p}^{nr}/{\goth
p}^\eta$
for some $\eta$. Let $\Cal E$ be an ordinary Anderson module over
$\overline{\n F_q[\th]/{\goth p}}$. Let us consider (3.1.1) for $\Cal E$, and let
$T_{\goth p}(\Cal E)$ be the Tate module of the etale part:
\medskip
$T_{\goth p}(\Cal E)=\underset{c \to\infty}\to{\invlim}\Cal E[{\goth p}^c]_{et}.$
\medskip
The dual Anderson module $\Cal E^t$ is defined in [L], [F], it is of rank $r$ and dimension $n$.
\medskip
Let $\Cal E_R$ be an Anderson module over $R$ such that its reduction to
$\overline{\n F_q[\th]/{\goth p}}$ is $\Cal E$ (a lift of $\Cal E$ on $R$).
\medskip
The function field analog of [K], Th. 2.1, (1) is the following:
\medskip
{\bf Theorem 3.4.1a.} The set of $\Cal E_R$ is in 1 -- 1 correspondence with the
set of maps
$$\Hom(T_{\goth p}(\Cal E)\otimes T_{\goth p}(\Cal E^t), \goth m),\eqno{(3.4.2)}$$

where Hom is of $\n F_q[T]/\goth p$-modules. 
\medskip
{\bf Notation.} For a fixed $\Cal E_R$ we denote this map by $q_{\Cal E_R}$.
\medskip
{\bf Idea of the proof. } First, we define the analog of the map $\vf _{\n A/R}$, [K], p. 151 for the present situation. Here it is $\vf _{\Cal E/R}: T_\goth p(\Cal E)\to \goth m^{\oplus n}$.
\medskip
Recall that $\eta$ satisfies $\goth m^\eta=0$. We choose $k$ such that $q^k\ge \eta$, and we consider formulas of multiplication by $\goth p^k$ for $\Cal E$:

$$\goth p^k (X)=\sum_{i=k}^\eta C_iX^{q^i}\eqno{(3.4.2a)}$$
where $X\in \p^{\oplus n}$ is a column vector and $C_i\in M_{n\times n}(\p)$. Condition $X\in \Cal E[\goth p^k]_{et} $ means that $\sum_{i=k}^\eta C_iX^{q^i}=0$.

Let $\tilde X\in R^{\oplus n}$ be a lift of  $X\in \Cal E[\goth p^k]_{et} $. Since for the first term $C_kX^{q^k}$ of 3.4.2a we have $q^k\ge N$, we get that $\vf _{\Cal E/R}(X):=\sum_{i=k}^\eta C_i\tilde X^{q^i}\in \goth m^{\oplus n}$ does not depend on the choice of $\tilde X$. $\square$
\medskip
Now we need the function field analog of [K], Th. 2.1, (4). Let $\Cal E_1$, $\Cal
E_2$ be ordinary Anderson modules over $\overline{\n F_q[\th]/{\goth p}}$,
$\alpha: \Cal E_1 \to \Cal E_2$ a map and $\Cal E_{1,R}$, $\Cal E_{2,R}$ lifts
of $\Cal E_1$, $\Cal E_2$ on $R$. We denote by

$$T_{\goth p}(\alpha): T_{\goth p}(\Cal E_1) \to T_{\goth p}(\Cal
E_2)\eqno{(3.4.3)}$$

$$T_{\goth p}(\alpha^t): T_{\goth p}(\Cal E_2^t) \to T_{\goth p}(\Cal
E_1^t)\eqno{(3.4.4)}$$
the maps obtained by functoriality.
\medskip
{\bf Lemma 3.4.5.} A map $\alpha_R: \Cal E_{1,R} \to \Cal E_{2,R}$
such that its
reduction is $\alpha$ exists iff for any $x\in T_{\goth p}(\Cal E_1)$, $y\in
T_{\goth p}(\Cal E_2^t)$ we have

$$q_{\Cal E_2,R}(T_{\goth p}(\alpha)(x), y) = q_{\Cal E_1,R}(x, T_{\goth
p}(\alpha^t)(y)).\eqno{(3.4.6)}$$
and moreover if this condition is satisfied then $\alpha_R$ is unique. $\square$
\medskip
{\bf Lemma 3.4.6a.} (Conjectural statement). To define $E_2$ over
$\n F_q[\th]_{\goth p}$
is the same as to define a concordant system of $(E_2)_\eta$ over
$R_\eta$ (the
concordance condition is clear).
\medskip
{\bf Remark.} Obviously $E_2$ defines a concordant system of $(\goth
E_2)_\eta$.
But is the inverse really true? Maybe non-trivial automorphisms of $(\goth
E_2)_\eta$ give obstacles?
\medskip
Now we return to the proof of Proposition 3.4. We fix $\eta$, we take $\Cal
E_1=\goth E_1$, $\Cal E_2=\goth E_2$. According (3.2), there exist bases
\medskip
$e_1^t, \dots,e_n^t$, $e_{n+1}, \dots,e_{r}$, $f_1^t, \dots,f_n^t$, $f_{n+1},
\dots,f_{r}$ of
\medskip
$T_{\goth p}(\goth E_1^t)$, $T_{\goth p}(\goth E_1)$, $T_{\goth p}(\goth
E_2^t)$, $T_{\goth p}(\goth E_2)$ respectively such that the maps $T_{\goth
p}(\delta)$, $T_{\goth p}(\delta^t)$ in these bases are the following:
\settabs 7 \columns
\medskip
$T_{\goth p}(\delta)(e_{n+1})=\goth pf_{n+1}$

\+ & ... & (maps of type 1)\cr

$T_{\goth p}(\delta)(e_{n+j})=\goth pf_{n+j}$,
\medskip
\medskip
$T_{\goth p}(\delta)(e_{n+j+1})=f_{n+j+1}$

\+ & ... & (maps of type 2)\cr

$T_{\goth p}(\delta)(e_{r})=f_{r}$,
\medskip
\medskip
$T_{\goth p}(\delta^t)(f_{1}^t)=\goth pe_{1}^t$

\+ & ... & (maps of type 3)\cr

$T_{\goth p}(\delta^t)(f_{i}^t)=\goth pe_{i}^t$,
\medskip
\medskip
$T_{\goth p}(\delta^t)(f_{i+1}^t)=e_{i+1}^t$

\+ &... & (maps of type 4)\cr

$T_{\goth p}(\delta^t)(f_{n}^t)=e_{n}^t$.
\medskip
\medskip
Now we apply formula (3.4.6) to these formulas. We consider 4 types of
$x$, $y$:

\medskip
{\bf Type 13.} $x$ of type 1, $y$ of type 3 ($\lambda\in[n+1,\cdots,n+j]$,
$\mu\in[1,\cdots,i])$:
\medskip
We get:
$$\goth p \cdot q_{\goth E_2,R}(f_\lambda, f_\mu^t)=\goth p \cdot q_{\goth
E_1,R}(e_{\lambda}, e_\mu^t).\eqno{(3.4.7)}$$

If $\goth m$ had no $\goth p$-torsion then we can divide the above equality by
$\goth p$ and to get
$$q_{\goth E_2,R}(f_\lambda, f_\mu^t)=q_{\goth E_1,R}(e_\lambda,
e_\mu^t),\eqno{(3.4.8)}$$

this means that $q_{\goth E_2,R}$ on these $f_\lambda$, $f_\mu^t$ is defined
uniquely.
\medskip
We think that in order to prove that we can really divide by $\goth p$, we must
consider not one fixed $\eta$, but all the values of them. The similar problem
exists for the next type:
\medskip
{\bf Type 14.} $x$ of type 1, $y$ of type 4 ($\lambda\in[n+1,\cdots,n+j]$,
$\mu\in[i+1,\cdots,n])$:
\medskip
We get: $$\goth p \cdot q_{\goth E_2,R}(f_{\lambda}, f_\mu^t)=q_{\goth
E_1,R}(e_{\lambda}, e_\mu^t).\eqno{(3.4.9)}$$

If $\goth m$ were $\goth p$-divisible and had the $\goth p$-torsion isomorphic
to $\n F_q[\th]/\goth p$ then OK: we have $\goth q^{(n-i)j}$ possibilities for $(\goth
E_2)_R$ as it should be.
\medskip
For other types of $x$, $y$ there is no such problem. Really:
\medskip
{\bf Type 23.} $x$ of type 2, $y$ of type 3 ($\lambda\in[n+j+1,\cdots,r]$,
$\mu\in[1,\cdots,i])$:
\medskip
We get: $$q_{\goth E_2,R}(f_{\lambda}, f_\mu^t)=\goth p \cdot q_{\goth
E_1,R}(e_{\lambda}, e_\mu^t).\eqno{(3.4.10)}$$

This means that $q_{\goth E_2,R}$ on these $f_\lambda$, $f_\mu^t$ is defined
uniquely;
\medskip
{\bf Type 24.} $x$ of type 2, $y$ of type 4 ($\lambda\in[n+j+1,\cdots,r]$,
$\mu\in[i+1,\cdots,n])$:
\medskip
We get: $$q_{\goth E_2,R}(f_{\lambda}, f_\mu^t)=q_{\goth E_1,R}(e_{\lambda},
e_\mu^t).\eqno{(3.4.11)}$$

This means that $q_{\goth E_2,R}$ on these $f_\lambda$, $f_\mu^t$ is defined
uniquely;
\medskip
{\bf 3.4.12.} We get that we have $\goth q^{(n-i)j}$ modules $E_2$,
this number
is equal to the quantity of $\Cal W$ satisfying (3.4.1).
\medskip
{\bf 3.4.13.} Now we need to prove that these $\Cal W$ really satisfy (3.4.1).
$\square$
\medskip
Now we can define the map $\beta_2:\n H(M)\to \Corr(\tilde X)$ from (1.1.4).
Idea of the definition: let $\tau_{\goth p}$ have the form $\left( \matrix
\goth p^A & 0 \\ 0 & \goth p^B  \endmatrix \right)$ where
$A=(a_1,\dots,a_{r-n})$, $B=(b_1,\dots,b_n)$, $\goth p^A =\diag(\goth
p^{a_1},\dots,\goth p^{a_{r-n}})$, $\goth p^b =\diag(\goth p^{b_1},\dots,\goth
p^{b_n})$. We denote by $T_{M,\goth p}=T_{M,\goth p}(A,B)$ the element
of Hecke
algebra $\n H(M)$ corresponding to the double coset $K_M\tau_{\goth p}K_M$.
Explicit formula for $S^{G_{fun}}_M$ shows that

$$S^{G_{fun}}_M(T_{G,\goth p})=\goth q^{-m_{A,B}}T_{M,\goth p} + \hbox{ other
terms,}\eqno{(3.5)}$$
where these other terms are linear combinations of $T_{M,\goth p}(A',B')$ for
pairs $(A',B')$ distinct from $(A,B)$. Coefficient $m_{A,B}\ge 0$ can
be easily
found explicitly; for $T_{M,\goth p}(A,B)=\Psi_i\Phi_j$ we have $m_{A,B}=ij$.
\medskip
Now we consider the reduction of the correspondence $\beta_1(T_{G,\goth p})$.
Let $\Gamma_{irr}$ be an irreducible component of its graph, $\phi:
E_1\to E_2$
a map of Anderson modular objects over $A/\goth p$ corresponding to
a point of
$\Gamma_{irr}$, and $\goth t\in \n H(M)$ its type. $\goth t$ depends only on
$\Gamma_{irr}$ but not on $\phi: E_1\to E_2$ because it is a discrete
invariant, so we can call it the type of $\Gamma_{irr}$.
\medskip
First, we denote by $\goth C(A,B)$ the correspondence on $\tilde X$ whose
graph is the sum of all the irreducible components of the graph of
reduction of
the correspondence $\beta_1(T_{G,\goth p})$ whose type is $T_{M,\goth p}(A,B)$
(really, for each $(A,B)$ there exists only one such component). By abuse of
notations we denote by $\goth C(\Psi_i\Phi_j)$ the $\goth C(A,B)$ where $A$,
$B$ are from 2.2.3 a,b. Finally, we define
$$\beta_2(T_{M,\goth p}(A,B))=\goth q^{m_{A,B}}\goth C(A,B),\eqno{(3.6)}$$
hence
$$\beta_2(\Psi_i\Phi_j)=\goth q^{ij}\goth C(\Psi_i\Phi_j).\eqno{(3.7)}$$

(3.5) and (3.6) show immediately that the function field analog of the diagram
(1.1.6) is commutative.
\medskip
{\bf Corollary 3.8.} $d_1^{ns}(\goth C(\Psi_i\Phi_j))=\goth q^{(n-i)j}$.
\medskip
{\bf Proof.} Follows immediately from 3.4.12. $\square$
\medskip
3.7 and 3.8 imply that
$$d_1^{ns}(\beta_2(\Psi_i\Phi_j))=\goth q^{nj}.\eqno{(3.9)}$$

{\bf Proposition 3.10.} $\beta_2$ is a ring homomorphism.
\medskip
{\bf Idea of the proof.} Let $(A,B)$, $(A',B')$ be 2 pairs of multiindices as
above and let $T_{M,\goth p}(A,B)\cdot T_{M,\goth p}(A',B')=\sum_i \kappa_i
T_{M,\goth p}(A_i,B_i)$ for some pairs $(A_i,B_i)$ and coefficients
$\kappa_i$.

\medskip
{\bf Lemma 3.10.1.} For all $i$ we have
$$\goth q^{m_{A,B}}d_1^{ns}(\goth C(A,B)) \cdot \goth
q^{m_{A',B'}}d_1^{ns}(\goth C(A',B'))=\kappa_i\goth
q^{m_{A_i,B_i}}d_1^{ns}(\goth C(A_i,B_i)).$$

{\bf Proof.} Explicit calculation. For a particular case corresponding to
$(A,B)=\Psi_i$, $(A',B')=\Phi_j$, $(A_i,B_i)=\Psi_i\Phi_j$ this follows from
the above results.

We have: $$T_{G,\goth p}(A,B)\cdot T_{G,\goth p}(A',B')=\sum_i \kappa_i
T_{G,\goth p}(A_i,B_i) + \hbox{ other terms}.\eqno{(3.10.2)}$$

Since the
reduction is a ring homomorphism, we see that:
\medskip
(a) (3.5) applied to the pairs $(A,B)$, $(A',B')$, $(A_i,B_i)$;
\medskip
(b) (3.10.2) and Lemma 3.10.1;
\medskip
(c) Commutativity of the the function field analog of the diagram (1.1.6)
\medskip
imply that
$$\beta_2(T_{M,\goth p}(A,B))\cdot \beta_2(T_{M,\goth p}(A',B'))+ \hbox{ other
terms}=$$ $$=\sum_i \kappa_i\beta_2(T_{M,\goth p}(A_i,B_i))+ \hbox{ other
terms.}$$

{\bf 3.10.3.} Now naive considerations show us that ``other terms'' in both
sides of the above equality are equal. Really, let us denote by $S\Gamma(A,B)$
the support of the graph $\Gamma(\beta_2(T_{M,\goth p}(A,B)))\subset \tilde X
\times \tilde X$, and analogically for the pairs $(A',B')$, $(A_i,B_i)$. We
have:
$$(t_1,t_2)\in S\Gamma(A,B) \iff \hbox{ there is a map } E_{t_1} \to E_{t_2}
\hbox{ of type } T_{M,\goth p}(A,B).$$

By definition of the product of
correspondences,
$$(t_1,t_3)\in \cup_i S\Gamma(A_i,B_i) \iff  \hbox{ there exists } t_2 \hbox{
such that } $$
$$(t_1,t_2)\in S\Gamma(A,B), \ \ (t_2,t_3)\in S\Gamma(A',B').$$

Since the type of the composition of maps of Anderson varieties is
concordant with the multiplication in $\n H(M)$, we get 3.10.3. $\square$
\medskip
\medskip
{\bf 3a. Conjectural form of Langlands correspondence. }\nopagebreak
\medskip
According Langlands, L-function $L(\Cal M,s)$ of an irreducible
submotive $\Cal
M$ of a Shimura variety is related with $L(\pi,\goth r,s)$,
where $\pi$ is an automorphic representation of $G(\n A_{\n Q})$ and $\goth r:$ $^L G \to GL(\goth W)$ a
finite-dimensional representation of $^L G$:
$$L(\Cal M,s)\sim L(\pi,\goth r,s).\eqno{(3a.0)}$$

Conjectural construction of $\goth r$ is given for example in [BR], Section 5.1, p. 550.
\medskip
Let us formulate an analog of this result for Anderson modular objects. For this case an analog of $G(\n A_{\n Q})$ is $G_{fun}(\n A_{\n F_q(\th)})$.
\medskip
{\bf Theorem 3a.1.} If an analog of (3a.0) is true for Anderson modular objects $X_{fun}(r,n)$ of any level
then the restriction of $\goth r$ to $\hat G_{fun}$ is the $n$-th skew power
representation of $GL_r$.
\medskip
This theorem follows from the below Theorem 3a.3.
\medskip
Let $\pi=\otimes \pi_{\goth l}$ be a representation of $G_{fun}(\n A_{\n F_q(\th)})$
corresponding (according Langlands) to an irreducible submotive of an Anderson modular object, and
$\theta_{\goth p}\in \ ^L G$ a Langlands element of $\pi_{\goth p}$
(we consider
the case of $\goth p$ such that $\pi_{\goth p}$ is non-ramified). Let $\alpha_i$, $i=1,...r$, be
eigenvalues of $\theta_{\goth p}$ and $a_i$ the eigenvalues of $T_{\goth p,i}$
(analogs of Fourier coefficients of an automorphic form for the classical
case). Standard formalism of Langlands elements for $GL_r$ in the
non-ramified case
together with (2.2.2) shows that
$$a_i=\goth q^{-i(i-1)/2}\sigma_i(\alpha_*).\eqno{(3a.2)}$$

We denote by $P'_{r,n}$ the characteristic polynomial of $\goth
r(\theta_{\goth
p})$, it belongs to $\n Z[a_1,...,a_r][T]$ where $T$ is an abstract variable.
The following theorem follows immediately from (2.2.2), (3a.2) (like in the
number case):
\medskip
{\bf Theorem 3a.3.} $P_{r,n}=P'_{r,n}$ (after identification of $T$ and $fr$,
$a_i$ and $T_{\goth p,i}$). $\square$
\medskip
\medskip
{\bf 4. Unitary Shimura varieties.}
\medskip
We consider abelian varieties with multiplication by an imaginary quadratic field (abbreviation: MIQF). Let $K$ be such field, $X_{num}$ the corresponding
Shimura variety parametrizing abelian $r$-folds with multiplication by $K$, of signature $(r-n,n)$. We shall call them unitary Shimura varieties. The corresponding reductive group over $\n Q$ is $G=G_{num}=GU(r-n,n)$. We
have $\dim X_{num}=(r-n)n$. Let $p$ be a prime
inert in $K$; we shall consider $p$-Hecke correspondences and the reduction at $p$.
\medskip
{\bf Theorem 4.1.} $M$ for this case is the same as in Theorem 2.1.
\medskip
{\bf Corollary 4.2.} Satake maps for this case coincide with the ones for the
functional case (formulas (2.2.4), (2.2.5)).
\medskip
Let $A_t$ be as in Subsection 1.2. $(A_t)_p$ is an $r$-dimensional
vector space
over $\n F_{p^2}$. Let $D=D_{unitary}$ be as in Subsection 1.2.
\medskip
{\bf Theorem 4.3.} $D_{unitary}$ is a vector space over $\n F_{p^2}$ of
dimension $max \ (r-n,n)$.
\medskip
{\bf Remark 4.4.} There exists a symmetry between $n$ and $r-n$. Nevertheless,
here the analogy between functional and unitary case apparently is not
complete.
\medskip
{\bf Theorem 4.5.} Analog of the Theorem 2.3.5 (i.e. formula 2.3.6)
holds for the
unitary case (dimension is taken over $\n F_{p^2}$).
\medskip
{\bf Theorem 4.6.} ([BR], Section 5.1, p. 550, example (b)). Restriction of
$\goth r$ on $\hat G \subset \ ^L G$ is the same as in Theorem 3a.1.
\medskip
We think that analogs of Theorems 2.3.7, 2.3.8, 2.3.10 also hold for this
case.
\medskip
{\bf 5. Conjectural values of Hodge numbers.}\nopagebreak
\medskip
There are conjectural formulas for values of Hodge numbers $h^{ij}$ (over the fields of multiplications coming from Hecke correspondences) of
irreducible submotives of Shimura varieties (see, for example, [BR], Section
4.3, p. 548). For example, for the case of Siegel modular varieties of genus
$g$ (their dimension is $d_g=g(g+1)/2$) and for a generic pure submotive of
weight $d_g$ they are the following:
\medskip
{\bf Theorem 5.1.} $h^{i,d_g-i}=$ \{the quantity of subsets $(j_1, \dots,
j_\alpha)$
of the set $1,2,...,g$ such that $j_1+ ... +j_\alpha=i$\}, where $\alpha$ is
arbitrary.
\medskip
For other types of submotives the formulas for $h^{ij}$ are similar but more long.
\medskip
For example, for the case of unitary Shimura variety of Section 4 and for the same type of
submotives the formula is the following:
\medskip
{\bf Theorem 5.2.} $h^{i,(r-n)n-i}=$ \{the quantity of subsets $(j_1,
... j_n)$ of
the set $(1,2,...,r)$ such that $j_1+ ... +j_n - (1+...+n)=i$\}.
\medskip
By analogy between functional and unitary case we can conjecture that the same
formula holds for the functional case.
\medskip
{\bf 6. Non-ordinary Drinfeld modules.}
\medskip
For further applications we shall state two problems and give their conjectural
answers. Let us restrict ourselves by the case $n=1$ and the correspondence
$T_{\goth p,1}$. These problems are related with the description of
intersection of two irreducible components of the graph of $T_{\goth p,1}$ in
characteristic $\goth p$. Let $E$ be a Drinfeld module such that its reduction
is a generic non-ordinary, i.e. $\dim_{\n F_{\goth p}}(D_{fun}(E))=2$ is the
least possible. Let $t$ be the point on $X_{fun}$ corresponding to $E$ and $t', t''
\in T_{\goth p,1}(t)$.
\medskip
{\bf Question 6.1. Formulate analogs of conjectures 2.3.5, 2.3.8 for
this $t$.}
\medskip
{\bf Conjectural answer.}
\medskip
{\bf (a)} $\tilde t'$ is (the only) closed point of
$\Psi_1(\tilde t)$ iff $\Ker(\alpha_{t,t'})\subset D_{fun}$.
\medskip
All such $\tilde t'$ also belong to
$\Phi_1(\tilde t)$.
\medskip
{\bf (b)} $\tilde t'=\tilde t''$ as closed points iff the $\n F_{\goth p}$-linear
spans coincide:
$$<\Ker(\alpha_{t,t'}), D_{fun}>=<\Ker(\alpha_{t,t''}), D_{fun}>.$$

Now we consider a more special situation. Let $r$ be even, $L$ a quadratic
extension of $\n F_q(\th)$ such that $\goth p$ inert in $L/\n F_q(\th)$, and $E$ a generic
Drinfeld
module with multiplication by $L$. In this case $D_{fun}(E)$ is a
1-dimensional
$\n F_{\goth p^2}$-vector space. Let $t,t',t''$ be as above.
\medskip
{\bf Question 6.2. Formulate analogs of conjectures 2.3.5, 2.3.8 for
this $t$.}

\medskip
{\bf Conjectural answer. (a)} is the same as above, and in {\bf (b)}
we consider
$\n F_{\goth p^2}$-linear spans:
$$\tilde t'=\tilde t''\iff <\Ker(\alpha_{t,t'}), D_{fun}>_{\n F_{\goth
p^2}}=<\Ker(\alpha_{t,t''}), D_{fun}>_{\n F_{\goth p^2}}.$$
\medskip
\newpage
{\bf References}
\medskip
[BR] Blasius D., Rogawski J.D., Zeta functions of Shimura varieties.
In: Motives. Proc. of Symp. in Pure Math., 1994, v. 55, part 2, p.
525 - 571
\medskip
[FCh] Faltings, Gerd; Chai, Ching-Li, Degeneration of abelian varieties.
Ergebnisse der Mathematik und ihrer Grenzgebiete (3) [Results in Mathematics
and Related Areas (3)], 22. Springer-Verlag, Berlin, 1990. xii+316 pp.
\medskip
[F] Faltings, Gerd, Group schemes with strict $\Cal O$-action. Mosc.
Math. J.  2
(2002),  no. 2, 249--279.
\medskip
[GL17] Grishkov, A.; Logachev, D. Lattice map for Anderson t-motives: first approach. J. of Number Theory. 2017, vol. 180, p. 373 -- 402. https://arxiv.org/pdf/1109.0679.pdf
\medskip
[GL18] Grishkov, A.; Logachev, D. Duality of Anderson t-motives having $ N\ne 0$.
https://arxiv.org/pdf/1812.11576.pdf
\medskip
[GL21] Grishkov, A.; Logachev, D. Anderson t-motives and abelian varieties with MIQF: results coming from an analogy. To appear in J. of algebra and its applications. https://arxiv.org/pdf/0907.4712.pdf
\medskip
[HJ] Hartl U.; Juschka A.-K. Pink's theory of Hodge structures and the Hodge conjecture over function fields. In: "t-motives: Hodge structures, transcendence and other motivic aspects", Editors G. B\"ockle, D. Goss, U. Hartl, M. Papanikolas, European Mathematical Society Congress Reports 2020, and https://arxiv.org/pdf/1607.01412.pdf
\medskip
[K] Katz, N. Serre-Tate local moduli.  Algebraic surfaces (Orsay, 1976--78),
pp. 138--202, Lecture Notes in Math., 868, Springer, Berlin-New York, 1981.
\medskip
[VL] Vincent Lafforgue, Param\`etres de Langlands et cohomologie des espaces de modules de $G$-chtoucas  arXiv:1209.5352
\medskip
[L] D. Logachev, Duality for abelian Anderson T-motives, and related problems.
http://arxiv.org/PS\_cache/arxiv/pdf/0711/0711.1928v4.pdf
\medskip
[P] Richard Pink, Hodge structures over function fields. Universit\"at Mannheim.
Preprint. September 17, 1997.
\medskip
[V] Yakov Varshavsky,  Moduli spaces of principal F-bundles. Sel. Math., New ser. 10 (2004) 131-166
\medskip
[W] T. Wedhorn, Congruence relations on some Shimura varieties,  Journ. f. d.
Reine und Angew. Math. 524 (2000) 43-71
\enddocument